\documentclass[draft]{amsart}
\usepackage{amssymb,graphicx,enumerate,amsthm}
\usepackage{multirow}

\newcommand{\RA}{\rightarrow}
\newcommand{\Sum}{\displaystyle\sum}

\renewcommand{\i}{\infty}

\newcommand{\beqs}{\begin{equation*}}
\newcommand{\eeqs}{\end{equation*}}
\numberwithin{equation}{section}
\newtheorem{theorem}{Theorem}[section]
\newtheorem{lemma}[theorem]{Lemma}

\newtheorem{definition}[theorem]{Definition}
\newtheorem{conjecture}[theorem]{Conjecture}

\begin{document}

\makeatletter
\def\imod#1{\allowbreak\mkern10mu({\operator@font mod}\,\,#1)}
\makeatother

\author{Alexander Berkovich}
   \address{Department of Mathematics, University of Florida, 358 Little Hall, Gainesville FL 32611, USA}
   \email{alexb@ufl.edu}

\author{Ali Kemal Uncu}
   \address{Department of Mathematics, University of Florida, 358 Little Hall, Gainesville FL 32611, USA}
   \email{akuncu@ufl.edu}

\title[\scalebox{.9}{A New Companion to Capparelli's Identities}]{A New Companion to Capparelli's Identities}
     
\begin{abstract} 
We discuss a new companion to Capparelli's identities. Capparelli's identities for $m=1,2$ state that the number of partitions of $n$ into distinct parts not congruent to $m, -m$ modulo $6$ is equal to the number of partitions of $n$ into distinct parts not equal to $m$, where the difference between parts is greater than or equal to $4$, unless consecutive parts are either both consecutive multiples of $3$ or add up to to a multiple of $6$. In this paper we show that the set of partitions of $n$ into distinct parts where the odd-indexed parts are not congruent to $m$ modulo $3$, the even-indexed parts are not congruent to $3-m$ modulo $3$, and $3l+1$ and $3l+2$ do not appear together as consecutive parts for any integer $l$ has the same number of elements as the above mentioned Capparelli's partitions of $n$. In this study we also extend the work of Alladi, Andrews and Gordon by providing a complete set of generating functions for the refined Capparelli partitions, and conjecture some combinatorial inequalities. \end{abstract}

  

\thanks{The first author is supported in part by The Simons Foundation Grant \#309829.}

\date{\today}
   
\maketitle
   
\section{Introduction and Notation}
\label{intro}

A partition $\pi$ is a finite, non-increasing sequence of positive integers $(\pi_1,\pi_2,\dots,\pi_k)$. The $\pi_i$ are called parts of the partition $\pi$, and $\pi_1+\pi_2+\dots+\pi_k$ is called the norm of $\pi$. We call $\pi$ "a partition of $n$'' if the norm of $\pi$ is $n$. Conventionally, we define empty sequence as the only partition of 0. Throughout this paper we assume that $a$ and $q$ are complex numbers where $|q|<1$, $L$ is a positive integer, and $m\in\{1,2\}$. 
We use the standard notations as in \cite{theoryofpartitions} and \cite{GasperRahman}:
	\begin{align}
	\nonumber
(a)_L&:=(a;q)_L = \prod_{n=0}^{L-1}(1-aq^n),\vspace*{.1cm}\\
	\nonumber
(a_1,a_2,\dots,a_k;q)_L &= (a_1;q)_L(a_2;q)_L\dots (a_k;q)_L,\vspace*{.1cm}\\
 \nonumber
  (a;q)_\i&:= \prod_{n=0}^{\infty}(1-aq^n).
   	\end{align}
	We define the $q$-binomial and $q$-trinomial coefficients, respectively, as
	\begin{align}
 \nonumber
  \displaystyle\genfrac{[}{]}{0pt}{}{k}{n}_q &:= \left\lbrace \begin{array}{ll}\frac{(q)_k}{(q)_n(q)_{k-n}}&\text{for }k\geq n \geq 0,\\
   0&\text{otherwise,}\end{array}\right.\\\intertext{and}
 \nonumber
  \displaystyle\genfrac{[}{]}{0pt}{}{k}{n,r}_q &:= \genfrac{[}{]}{0pt}{}{k}{n}_q\genfrac{[}{]}{0pt}{}{k-n}{r}_q= \left\lbrace \begin{array}{ll}\frac{(q)_k}{(q)_n(q)_r(q)_{k-n-r}}&\text{for }k\geq n+r\geq n\geq 0,\\
   0&\text{otherwise.}\end{array}\right.
   \end{align}

Let $C_m(n)$ be the number of partitions of $n$ into distinct parts where no part is congruent to $\pm m$ modulo $6$. Define $D_m(n)$ to be the number of partitions of $n$ into parts, not equal to $m$, where the minimal difference between consecutive parts is 2, in fact, the difference between consecutive parts is greater than or equal to $4$ unless these parts either are both consecutive multiples of 3 (yielding a difference of $3$) or add up to a multiple of 6 (possibly giving a difference of $2$).
   
In 1988, S.~Capparelli stated two conjectures for $C_m$ and $D_m$ in his thesis \cite{CapparelliThesis}. The first one was later proven by G.~E.~Andrews \cite{AndrewsCapparelli} in 1992 during the Centenary Conference in Honor of Hans Rademacher. Two years later Lie theoretic proof were supplied by Tamba and Xie \cite{Tamba} and by Capparelli \cite{CapparelliProof}. The first of Capparelli's conjectures was stated and proven in the form of Theorem~\ref{capparelliconj}.
\begin{theorem}[Andrews 1992]\label{capparelliconj} For any non-negative integer $n$, \begin{equation}\nonumber
C_1(n) = D_1(n).\end{equation}
\end{theorem}

As an example let $n=19$; and find that $C_1(19)=D_1(19)=10$. Table~\ref{table1} gives the partitions associated with $n=19$ applied in Theorem~\ref{capparelliconj}.
\begin{center}\begin{table}[htb]
\caption{$C_1$ and $D_1$ numbers and respective partitions for Theorem~\ref{capparelliconj} with $n=19$}\begin{tabular}{cc}
$C_1(19)=10$ : & $\begin{array}{c}
(16,3),\ (15,4),\ (14,3,2),\ (12,4,3),\ (10,9), \vspace*{.1cm}\\
(10,6,3),\ (10,4,3,2),\ (9,8,2),\ (9,6,4),\ (8,6,3,2)\vspace*{.1cm}
\end{array}$ \\ 
\\ [-1.5ex]
$D_1(19)=10$ : & $\begin{array}{c}\vspace*{.1cm}(19),\ (17,2),\ (16,3),\ (15,4),\ (14,5),\\\vspace*{.1cm} (13,6),\ (13,4,2),\ (12,7),\ (11,6,2),\ (10,6,3) \vspace*{.1cm}\end{array}$ 
\end{tabular}\label{table1} \end{table}\vspace*{-.5cm}\end{center}

One year after Capparelli's proof, Alladi, Andrews, and Gordon improved on Theorem~\ref{capparelliconj} in \cite{refinement}. They gave a refinement of these identities by introducing restrictions on the number of occurrences of parts belonging to certain congruence classes. A detailed discussion of their generating functions is presented in Section~\ref{Further Results}. In particular, they stated and proved the following extension of Theorem~\ref{capparelliconj}.

\begin{theorem}[Capparelli 1994; Alladi, Andrews, Gordon 1995]\label{fulcapparelli} For any non-negative integer $n$ and $m\in\{1,2\}$,\begin{equation}\nonumber
C_m(n) = D_m(n).\end{equation}
\end{theorem}

We can illustrate Theorem~\ref{fulcapparelli} by taking $m=2$, $n=19$ and listing the associated partitions in Table~\ref{table2}.
\begin{center}\begin{table}[htb]\caption{$C_2$ and $D_2$ numbers and respective partitions for Theorem~\ref{fulcapparelli} with $n=19$}\begin{tabular}{cc}
$C_2(19)=12$ :& $\begin{array}{c}
(19),\ (18,1),\ (15,3,1),\ (13,6),\ (13,5,1),\ (12,7), \vspace*{.1cm}\\
(12,6,1),\ (11,7,1),\ (11,5,3),\ (9,7,3),\ (9,6,3,1),\ (7,6,5,1) \vspace*{.1cm}
\end{array}$ \\ 
 \\ [-1.5ex]
$D_2(19)=12$ :& $\begin{array}{c}(19),\ (18,1),\ (16,3),\ (15,4),\ (14,5),\ (13,6), \vspace*{.1cm}\\ (13,5,1),\ (12,7),\ (12,6,1),\ (11,7,1),\ (10,8,1),\ (10,6,3) \vspace*{.1cm} \end{array}$  
\end{tabular}\label{table2}\end{table}\vspace*{-.5cm} \end{center}

The reader is referred to articles \cite{bringmann} by Bringmann and Mahlburg, \cite{Dousse} by Dousse, and \cite{sills} by Sills for recent research on Capparelli's identities.

Let $\lfloor x\rfloor$ denote the greatest integer less than or equal to $x$. Let $A_m(n)$ be the number of partitions $\pi=(\pi_1,\pi_2,
\dots)$ of $n$ such that
\begin{enumerate}[\itshape i.]
\item $\pi_{2i+r} \not\equiv 3-m+(-1)^m r\pmod 3$, and
\item $\pi_{2i+r} - \pi_{2i+1+r}>\lfloor m/2\rfloor + (-1)^{m-1}r$ for $1\leq 2i+r $
\end{enumerate}
where $r\in\{0,1\}$. We remark that the second condition of the definition of $A_m(n)$ can be replaced with the condition that \textit{all parts are distinct and $3l+1$ and $3l+2$ do not appear together as consecutive parts for any integer $l$}.

Theorem~\ref{newcapparellicomp} is a new companion to Theorem~\ref{fulcapparelli},
\begin{theorem}\label{newcapparellicomp} Let $n$ be a non-negative integer and $m\in\{1,2\}$. Then, \begin{equation}\nonumber
A_m(n) = C_m(n).\end{equation}
\end{theorem} 

We remark that one can obtain this new companion, Theorem~\ref{newcapparellicomp} as a non-trivial corollary of Boulet's results in \cite{Boulet}. The details will be given elsewhere.

Comparing Table~\ref{table3} with Table~\ref{table1} and Table~\ref{table2} leads us to the $n=19$ case for Theorem \ref{newcapparellicomp}. 
\begin{center}\begin{table}[htb]\caption{$A_1$ and $A_2$ numbers and respective partitions for Theorem~\ref{newcapparellicomp} with $n=19$}
\begin{tabular}{cc}
$A_1(19)=10$ : & $\begin{array}{c}
(18,1),\ (15,4),\ (14,3,2),\ (12,7),\ (12,4,3), \vspace*{.1cm} \\
(11,6,2),\ (11,4,3,1),\ (9,7,3),\ (9,6,3,1),\ (8,6,5). \vspace*{.1cm}
\end{array}$ \\ 
 \\ [-1.5ex]
$A_2(19)=12$ :& $\begin{array}{c}(19),\ (16,3),\ (15,3,1),\ (13,6),\ (13,5,1),\ (12,6,1),\\ \vspace*{.1cm} (10,9),\ (10,8,1),\ (10,6,3),\ (10,5,4),\ (9,6,4),\ (9,5,3,2).  \vspace*{.1cm}\end{array}$ 
\end{tabular}\label{table3}
\end{table}\vspace*{-.5cm} \end{center}

In this paper we discuss a refined version of Theorem~\ref{newcapparellicomp}, which is our main combinatorial result Theorem~\ref{refined}, and then obtain Theorem~\ref{newcapparellicomp} as a consequence. We give the necessary definitions and discuss recurrences in Section~\ref{GenFuncs}. Section~\ref{SecProofMain} is reserved for the proof of the Theorem~\ref{genfuncsdef}. We provide combinatorial interpretation of \eqref{def1} in Section~\ref{CapparelliCompSec} and relate this result with our refinement of Capparelli's companion, Theorem~\ref{newcapparellicomp}. In Section~\ref{Further Results} we discuss three new polynomials which are generating functions for number of partitions with Capparelli-type difference conditions subject to the various bounds on the largest parts of partitions. This gives the full list of generating functions and extends the study of Alladi, Andrews and Gordon \cite{refinement}. We conclude this section with some $q$-theoretic and combinatorial conjectures.

\section{Generating Functions and Recurrences}\label{GenFuncs}

In this section we assume that $m\in\{1,2\}$; $N,\ i$, and $j$ are non-negative integers; and $N\geq i,\ j $. Let $x$ be a real number and define $\lceil x\rceil$ to be the least integer greater or equal to $x$; $\lfloor x\rfloor$ to be the greatest integer less than or equal to $x$; and $\{\{x\}\}$ to be the fractional part of $x$.
\begin{definition}\label{Pgen}
Define $P_{m,N}(i,j,q)$ to be the generating function for the number of partitions $\pi = (\pi_{1},\pi_{2},\dots)$ such that
\begin{enumerate}[\itshape i.]
\item $\pi_{2l+r} \not\equiv 3-m+(-1)^m r\pmod 3$, 
\item $\pi_{2l+r} - \pi_{2l+1+r}>\lfloor m/2\rfloor + (-1)^{m-1}r$ for $1\leq 2l+r$.
\item $\pi_{1}$, is $\leq 3\lceil N/2 \rceil  - 2m\{\{N/2\}\}$,
\item $i$ is the number of parts congruent to $2\pmod 3$, and
\item $j$ is the number of parts congruent to $1\pmod 3$
\end{enumerate} 
where $r\in\{0,1\}$.
\end{definition}

Let $A_{m, N}(n,i,j)$ be the number of partitions of $n$ satisfying the conditions i--iv of Definition~\ref{Pgen}. From Definition~\ref{Pgen}, it is easy to see that
\begin{align}
\nonumber
\lim_{N\RA\i} \Sum_{i,j=0}^\i P_{m,N}(i,j,q) &= \lim_{N\RA\i} \Sum_{n,i,j=0}^\i A_{m,N}(n,i,j)q^n =  \Sum_{n=0}^\i A_m(n) q^n,\\
\intertext{where} 
\label{A_mN}
A_m(n) &= \lim_{N\rightarrow\infty}\Sum_{i,j=0}^{\infty} A_{m,N}(n,i,j).
\end{align}
The series in \eqref{A_mN} is a finite sum, because $A_{m,N}(n,i,j)=0$ for all $i$ and $j\geq n$.
The explicit expression of these generating functions $P_{m,N}(i,j,q)$ are presented in Theorem~\ref{genfuncsdef}. This is the main theorem of this paper.

\begin{theorem}\label{genfuncsdef}
For non-negative integers $N,i,j$ where $N\geq i,j$, and $m=1,2$, we have
\begin{align}
\label{def1}P_{m,2N} (i,j,q) &= q^{\omega(m,i,j)} \genfrac{[}{]}{0pt}{}{N}{i,\ j}_{q^6} (-q^3;q^3)_{N-i-j},\\
\label{def2}P_{m,2N+1} (i,j,q) &= q^{\omega(m,i,j)} \frac{1-q^{3(N+1+i+(-1)^m j)}}{1-q^{6(N+1)}} \genfrac{[}{]}{0pt}{}{N+1}{i,\ j}_{q^6} (-q^3;q^3)_{N+1-i-j},
\end{align}
where $\omega(m,i,j) := (3i+(-1)^m m)i + (3j+(-1)^{m+1}m)j$.
\end{theorem}

In order to prove Theorem~\ref{genfuncsdef} we need the following recurrences:

\begin{lemma}\label{recurrence}
For $N,i,\ $and$\ j$ defined as before,
\begin{align}
\label{rec1}P_{1,2N+1} (i,j,q) &= P_{1,2N} (i,j,q) + \chi(i>1) q^{3N+2}P_{2,2N}(i-1,j,q),\\
\label{rec2}P_{1,2N+2} (i,j,q) &= P_{1,2N+1}(i,j,q) + q^{3(N+1)}P_{2,2N+1}(i,j,q),\\
\label{rec3}P_{2,2N+1} (i,j,q) &= P_{2,2N}(i,j,q) + \chi(j>1) q^{3N+1} P_{1,2N} (i,j-1,q),\\\intertext{and}
\label{rec4}P_{2,2N+2} (i,j,q) &= P_{2,2N+1} + q^{3(N+1)}P_{1,2N+1}(i,j,q),
\end{align}
\end{lemma}
where \begin{equation}\nonumber\chi(\text{statement})= \left\lbrace \begin{array}{ll}
1 & \text{if statement is true},\\
0 &\text{otherwise.}
\end{array} \right.\end{equation}

Lemma~\ref{recurrence} along with the initial conditions $P_{m,0}(i,j,q) = \delta_{i,0}\delta_{j,0}$ for $m=1,\ 2$ uniquely specifies these generating functions. Here the Kronecker delta function $\delta_{i,j} = 1$ if $i=j$, and $0$, otherwise. Similar to Definition~\ref{Pgen} we define generating functions for the number of partitions for a particular refinement of Capparelli-type congruence conditions.

\begin{definition}\label{Qgen}
For $m=1, 2$, let $Q_{m,N}(i,j,q)$ be the generating function for the number of partitions into distinct parts where
\begin{enumerate}[i.]
\item no part is congruent to $\pm m\pmod6$,
\item there are exactly $i$ parts $\equiv m+(-1)^{m+1} \pmod6$ and these parts are all $\leq6N-(3+m)$, 
\item there are exactly $j$ parts $\equiv 3+m \pmod6$ and these parts are all $\leq6(N-i)-(m+(-1)^{m+1})$, 
\item all parts that are $0\pmod3$ are bounded by $3(N-i-j)$.
\end{enumerate} 
\end{definition}

It is clear from Definition~\ref{Qgen} that the bounds on the parts depend on the congruence classes modulo 6. We proceed by formulating the main generalization of the companion result to Capparelli's identities. Let $C_{m,N}(n,i,j)$ be the number of partitions of $n$ satisfying the conditions i--iv of Definition~\ref{Qgen}; explicitly, \begin{equation}\nonumber
Q_{m,N}(i,j,q) =\Sum_{n=0}^{\infty} C_{m,N}(n,i,j) q^n.
\end{equation}
The refinement of Theorem~\ref{newcapparellicomp} is the following theorem:
\begin{theorem}\label{refined}
For $N,n,i,j\in\mathbb{Z}_{\geq0}$ and $m\in\{1,2\}$, 
\begin{equation}\nonumber
A_{m,2N}(n,i,j) = C_{m,N}(n,i,j).\end{equation}
\end{theorem}
Theorem~\ref{refined} implies Theorem~\ref{newcapparellicomp} by summing over $i$ and $j$ and letting $N$ tend to infinity.

\section{Proof of Theorem~\ref{genfuncsdef}}\label{SecProofMain}

In this section we prove  Theorem~\ref{genfuncsdef} which is one of the main results of this paper. We begin by proving Lemma~\ref{recurrence}. 

\begin{proof} Let $m=1$, and $N,\ i,$ and $j$ be non-negative integers satisfying $N\geq i,\ j$. The first recursion, \eqref{rec1}, comes from elementary observations. Let $\pi  = (\pi_1,\pi_2,\dots,\pi_k)$ be a partition satisfying the conditions in Definition~\ref{Pgen} with $m=1$, and $N\mapsto 2N+1$. If $\pi_1$ less than $3N+2$ then $\pi$ must also satisfy the conditions for $P_{1,2N}(i, j,q)$ because the only difference between $P_{m,2N+1}(i,j,q)$ and $P_{m,2N}(i,j,q)$ is in the bounds on the largest parts.  If $\pi_1 = 3N+2$ (which implicitly requires $i>0$) we can extract this part from $\pi$ and get a new partition. The leftover partition $\pi'~= (\pi_2,\pi_3,\dots,\pi_k)~=~(\pi'_1, \pi'_2,\dots,\pi'_{k-1})$ has one less count of $2$ modulo $3$ parts ($i\mapsto i-1$) and the largest part of $\pi'$, $\pi'_1$, is bounded by $3N$. Lastly the congruence conditions ii and iii in Definition~\ref{Pgen} for $P_{2,2N}(i-1,j,q)$ are satisfied by $\pi'$ as the extraction of the largest part from $\pi$ alters the parities of the indices of parts. Hence, we get the recurrence $P_{1,2N+1} (i,j,q) = P_{1,2N} (i,j,q) + \chi(i>0) q^{3N+2}P_{2,2N}(i-1,j,q)$.

The recurrences \eqref{rec2}, \eqref{rec3}, and \eqref{rec4} can similarly be established by examining the partitions satisfying the conditions for their respective definitions.
\end{proof}

In order to prove Theorem~\ref{genfuncsdef} we need to show that both sides of the equations \eqref{def1} and \eqref{def2} satisfy the same recurrences \eqref{rec1}--\eqref{rec4} with the same initial conditions. The recurrences of the left-hand side of the equations in Theorem~\ref{genfuncsdef} are handled in Lemma~\ref{recurrence}. Next, we show that the right-hand side of the equations in Theorem~\ref{genfuncsdef} satisfy the recurrences of Lemma~\ref{recurrence}.

\begin{proof}We will start with the right-hand side of \eqref{def1}. Let $m=1$, $N,\ i$, and $j$ be non-negative integers satisfying $N\geq i+j$. Then,
\begin{align}
\label{RHSofREC}q^{(3i-1)i + (3j+1)j}&\genfrac{[}{]}{0pt}{}{N}{i,\ j}_{q^6}(-q^3;q^3)_{N-i-j}\\[0.3cm]
\nonumber
&= q^{(3i-1)i + (3j+1)j}\genfrac{[}{]}{0pt}{}{N}{i,\ j}_{q^6} (-q^3;q^3)_{N-i-j} \frac{1-q^{3(N+i-j)}+q^{3(N+i-j)}+q^{6N}}{1-q^{6N}}\\[0.3cm]
&=\label{explicitreccurrence} \left(q^{(3i-1)i + (3j+1)j}\frac{1-q^{3(N+i-j)}}{1-q^{6N}}\genfrac{[}{]}{0pt}{}{N}{i,\ j}_{q^6} (-q^3;q^3)_{N-i-j}\right)\\[0.3cm] \nonumber&\hspace{.5cm} +
q^{3N}\left(q^{(3i+2)i + (3j-2)j}\frac{1-q^{3(N-i+j)}}{1-q^{6N}}\genfrac{[}{]}{0pt}{}{N}{i,\ j}_{q^6} (-q^3;q^3)_{N-i-j}\right).
\intertext{Comparison between \eqref{RHSofREC} and \eqref{explicitreccurrence} shows that \eqref{RHSofREC} satisfies the same recursion relation as $P_{1,2N}(i,j,q)$ given in \eqref{rec2}. Similarly the right-hand side of \eqref{def2} with $m=1$ satisfies the recurrence \eqref{rec1}. Here the $i=0$ case is obvious as the recurrences reduce down to $q$-binomial recurrences. Suppose that $i\geq 1$:
}
q^{(3i-1)i + (3j+1)j}&\frac{1-q^{3(N+1+i-j)}}{1-q^{6N}}\genfrac{[}{]}{0pt}{}{N+1}{i,\ j}_{q^6} (-q^3;q^3)_{N+1-i-j}\\[0.3cm] 
\nonumber
&= q^{(3i-1)i + (3j+1)j}\genfrac{[}{]}{0pt}{}{N+1}{i,\ j}_{q^6} (-q^3;q^3)_{N+1-i-j}\ \\[0.3cm] \nonumber &\hspace*{.5cm}\times \left(\frac{1-q^{3(N+1-i-j)}+q^{3(N+1-i-j)}+q^{3(N+1+i-j)}}{1-q^{6(N+1)}} \right)\\[0.3cm]
&= \left(q^{(3i-1)i + (3j+1)j} \genfrac{[}{]}{0pt}{}{N}{i,\ j}_{q^6} (-q^3;q^3)_{N-i-j} \right)\\[0.3cm] \nonumber&\hspace*{.5cm} + q^{3N+2} \left(q^{(3i-1)(i-1) + (3j-2)j} \genfrac{[}{]}{0pt}{}{N+1}{i-1,\ j}_{q^6} (-q^3;q^3)_{N+1-i-j} \right).
\end{align}
Analogous proofs can be easily given for the right-hand sides of \eqref{def1} and \eqref{def2} with $m=2$ in Theorem~\ref{genfuncsdef}. Picking $N=i=j=0$ in the right-hand side of the \eqref{def1} we see that these functions have the same initial conditions as $P_{m,N}(i,j,q)$, which finishes the proof of Theorem~\ref{genfuncsdef}.
\end{proof}

\section{Companions to Capparelli's Identities}\label{CapparelliCompSec}

We suppose that $m\in\{1,2\}$, $N,\ i$, and $j$ are non-negative integers where $N\geq i,\ j $. The proof of Theorem~\ref{refined} follows from showing that $P_{m,2N}(i,j,q)$, and $Q_{m,N}(i,j,q)$ are equal. We focus our attention on the product representation \eqref{def1} of the generating functions $P_{m,2N}(i,j,q)$.
For $m=1$, he expression
\begin{equation}\nonumber
q^{(3i-1)i + (3j+1)j} \genfrac{[}{]}{0pt}{}{N}{i,\ j}_{q^6} (-q^3;q^3)_{N-i-j},
\end{equation}
can be rewritten in terms of $q$-binomial coefficients as
\begin{equation}\nonumber
({q^{6}})^{{i+1 \choose 2}} \genfrac{[}{]}{0pt}{}{N}{i}_{q^6} (q^{6})^{{j+1 \choose 2}} \genfrac{[}{]}{0pt}{}{N-i}{j}_{q^6} (-q^3;q^3)_{N-i-j}\ q^{-4i-2j}.
\end{equation}
The factor
\begin{equation}\label{qbin_bounded_pair}
({q^{6}})^{{i+1 \choose 2}} \genfrac{[}{]}{0pt}{}{N}{i}_{q^6}
\end{equation}
is the generating function for the number of partitions into $i$ distinct multiples of 6 less than or equal to $ 6N$. Multiplying \eqref{qbin_bounded_pair} with the term $q^{-4i}$ can be interpreted as taking off 4 from each and every one of the $i$ parts. Consequently,
\begin{equation}\nonumber
q^{-4i}({q^{6}})^{{i+1 \choose 2}} \genfrac{[}{]}{0pt}{}{N}{i}_{q^6}
\end{equation}
is the generating function for the number of partitions into $i$ distinct parts less than or equal to $6N-4$ where every part is congruent to $2$ modulo $6$.
Similarly,
\begin{equation}\nonumber
q^{-2j}(q^{6})^{{j+1 \choose 2}} \genfrac{[}{]}{0pt}{}{N-i}{j}_{q^6}
\end{equation}
is the generating function for the number of partitions into $j$ distinct parts less than or equal to $6(N-i)-2$ where every part is congruent to $4$ modulo $6$.
Finally, we also know that $(-q^3;q^3)_{N-i-j}$ is the generating function for number of partitions into distinct parts less than or equal to $3(N-i-j)$. Discussion of the $m=2$ case can be given along the similar lines. Therefore we get the proof of Theorem~\ref{refined} as follows:

\begin{proof} The above construction shows that \begin{equation}\nonumber
P_{1,2N}(i,j,q)= q^{(3i-1)i + (3j+1)j} \genfrac{[}{]}{0pt}{}{N}{i,\ j}_{q^6} (-q^3;q^3)_{N-i-j} = Q_{1,N}(i,j,q),
\end{equation}
thus giving us the refined companion to Capparelli's identity. $P_{2,2N}(i,j,q) = Q_{2,N}(i,j,q)$ can be shown in the same manner.
\end{proof}

\section{Further Observations}\label{Further Results}

In \cite{refinement}, Alladi, Andrews, and Gordon refined Capparelli's identity using $G_{1,3N-2}(a,b,q)$, which is to be defined later. In this section we discuss the dominance and subordinance properties between these generating functions compared to the previously defined generating functions $P_{m,N}(i,j,q)$ of Definition~\ref{Pgen}. This approach will be used to conjecture combinatorial inequalities between sets of partitions for any given positive integer $n$.

For $m\in\{1,2\}$, let $G_{m,N} (a,b,q)$ be the generating function for number of partitions where \begin{enumerate}[\itshape i.] 
\item parts are not equal to $m$, 
\item the difference between consecutive parts is greater or equal than $4$ unless either consecutive parts are both consecutive multiples of $3$ or add to a multiple of 6, 
\item largest part is less than or equal to $N$, \item the exponent of $a$ counts the number of parts congruent to $2$ modulo $3$, and 
\item the exponent of $b$ counts the number of parts congruent to $1$ modulo $3$.
\end{enumerate}
For a positive integer $N$, it is easy to see that these generating functions satisfy the recursion relations:
\begin{align}
\label{AAGrec1}G_{m,3N+1} &= G_{m,3N} + bq^{3N+1}G_{m,3(N-1)} + abq^{6N}G_{m,3(N-2)+1},\\
\label{AAGrec2}G_{m,3N} &= G_{m,3(N-1)+2} + q^{3N}G_{m,3(N-1)},\\
\label{AAGrec3}G_{m,3N+2} &= G_{m,3N+1} + aq^{3N+2}G_{3(N-1)+1}.
\end{align}
Moreover, the first four initial conditions 
\begin{equation}
\nonumber\begin{array}{cccl} G_{m,-2}=\delta_{1,m}, & G_{m,-1}=1, & G_{m,0}=1, &\text{and } G_{m,1} = 1 + \delta_{2,m}bq,
\end{array} \end{equation}
define both generating function sequences uniquely.

Combining \eqref{AAGrec1}--\eqref{AAGrec3} one gets the equality
\begin{equation}
\nonumber q^{3N}(1+bq)G_{m,3(N-1)} = G_{m,3N+1} - G_{m,3(N-1)+1}  - aq^{3N-1}(1+bq^{3N+1})G_{m,3(N-2)+1}.
\end{equation}
This result combined with the analogous result of Alladi, Andrews and Gordon (\cite{refinement}, (4.1)) gives us the recurrence
\begin{equation}
\label{afterAAg4.1} G_{m,3(N-1)} =\frac{ G_{m,3N+1} + bqG_{m,3(N-1)+1}  - abq^{3N}(1-q^{3N})G_{m,3(N-2)+1}}{(1+bq)}.
\end{equation}

The definition of $G_{m,3(N-1)}$ implies that it is a polynomial in $a$, $b$, and $q$. Therefore the numerator in \eqref{afterAAg4.1} is divisible by $1+bq$. Furthermore, with the choice of variables $m=1$, $a=1/t$ and $b=t$ in (\ref{AAGrec1} -- \ref{AAGrec3}) one gets the recurrences of Andrews in \cite{AndrewsCapparelli}; namely (\cite{AndrewsCapparelli}, (4.2)--(4.9)).

One can show that for $m=1,2$, $G_{m,3N}$ and $G_{m,3N-2}$ have explicit polynomial representations. Hence, \eqref{AAGrec3} provides a polynomial representation for $G_{m,3N-1}$. In their paper \cite{refinement}, Alladi, Andrews and Gordon found the polynomial representation for the generating function
\begin{equation}\label{polyAAG}G_{1,3N-2}(a,b,q) = \sum_{l=0}^{\lfloor N/2 \rfloor} (q^3)^{N-2l \choose  2} \genfrac{[}{]}{0pt}{}{N}{2l}_{q^3}(-aq^2,-bq^4;q^6)_l.\end{equation}
In fact, we discovered new formulas for all three remaining generating functions. In particular,
\begin{equation} G_{2,3N-2}(a,b,q) =\sum_{l=0}^{\lfloor N/2-1 \rfloor} (q^3)^{N-2l-2 \choose  2} \genfrac{[}{]}{0pt}{}{N-1}{2l+1}_{q^3}(-aq^5;q^6)_l(-bq;q^6)_{l+1}.\end{equation}
Moreover, defining polynomials 
\begin{align}
S(a,b,q,N) &= \sum_{l=0}^{\lfloor N/2 \rfloor} (q^3)^{N-2l \choose  2} \genfrac{[}{]}{0pt}{}{N+1}{2l+1}_{q^3}(-aq^2,-bq^4;q^6)_{l},\\
T(a,b,q,N) &= \sum_{l=0}^{\lfloor N/2 \rfloor} (q^3)^{N-2l \choose  2} \genfrac{[}{]}{0pt}{}{N+1}{2l+1}_{q^3}(-aq^5,-bq;q^6)_{l},
\end{align}
we have \begin{align}
G_{1,3N}(a,b,q) &= S(a,b,q,N) + aq^{3N-1} S(a,b,q,N-1),\\
G_{2,3N}(a,b,q) &= T(a,b,q,N) + bq^{3N-2} T(a,b,q,N-1).
\end{align}
Obviously, these new representations for $G_{m,3N}$ with $m=1,2$ are polynomials. 

Now, we define $\Psi_{m,N}(a,b,q)$ to be the generating function for the number of partitions counted by $P_{m,N}(i,j,q)$, where the exponent of $a$ counts the number of parts congruent to $2$ modulo $3$, and the exponent of $b$ counts the number of parts congruent to $1$ modulo $3$. This definition can be written as a single-fold sum as \begin{equation}\nonumber\Psi_{m,N}(a,b,q) = \sum_{i,j,=0}^{\infty} P_{m,N}(i,j,q)a^ib^j.\end{equation}
We have the following relations, \begin{equation}\nonumber\lim_{N\rightarrow\infty}\Psi_{m,N}(a,b,q) = \lim_{N\rightarrow\infty}G_{m,N}(a,b,q)=(-aq^{3-(-1)^{m+1}m},-bq^{3-(-1)^m m};q^6)_\infty(-q^3;q^3)_\infty.\end{equation}

There are restrictions on the largest part of a partition in both definitions of the generating functions $G_{m,N}(a,b,q)$ and $\Psi_{m,N}(a,b,q)$. These restrictions ensure that they are both polynomials of variables $a$, $b$, and $q$. The bounds on the largest parts become equal for different choices of $N$, as in the example of $\Psi_{m,2N}(a,b,q)$ and $G_{m,3N}(a,b,q)$. Neither the degrees nor the coefficients of these polynomial are identical. Moreover, the coefficients in the power series expansion of the difference of these polynomials are of the same sign, which leads into the discussion of dominance.

Given two series $\sum_{i,j,k} c_{i,j,k}a^ib^jq^k$ and $\sum_{i,j,k} d_{i,j,k}a^ib^jq^k$, we define $\sum c_{i,j,k}a^ib^jq^k \succeq \sum d_{i,j,k}a^ib^jq^k$ if $c_{i,j,k}\geq d_{i,j,k}$ for all $i,j$ and $k$. We call $\sum c_{i,j,k}a^ib^jq^k$ dominant over $\sum d_{i,j,k}a^ib^jq^k$ (or $\sum d_{i,j,k}a^ib^jq^k$ to be subordinate). With this definition we can write the Conjecture~\ref{conjecture}:

\begin{conjecture}\label{conjecture} For $m=1,\ 2$, \begin{align}\label{conj}\Psi_{m,2N}(a,b,q) &\succeq G_{m,3N}(a,b,q),\\
\Psi_{m,2N+1}(a,b,q) &\succeq G_{m,3N+3-m}(a,b,q).\label{conjj}
\end{align}
 For $m=1,\ 2$ and $N\geq 1$, these dominance relations can be written more explicitly as
\begin{align}
\nonumber\Psi_{m,2N}(a,b,q) - G_{m,3N}(a,b,q) &= a^{\delta_{2,m}}b^{\delta_{1,m}} q^{3N+m}+\dots,\\
\nonumber\Psi_{m,2N+1}(a,b,q) - G_{m,3N+3-m}(a,b,q) &= a^{\delta_{1,m}}b^{\delta_{2,m}} q^{(3(N+2)-1)\delta_{1,m} + (3(N+1)-2)\delta_{2,m}}+\dots,
\end{align} where there are only non-negative higher degree terms under the ellipsis.
\end{conjecture}

These observed dominance properties, \eqref{conj} and \eqref{conjj}, imply various combinatorial inequalities. One such inequality is given as the following statement:

\noindent
\textit{For  $m=1,\ 2$, the number of partitions of $n$ with parts less than or equal to $3N$ satisfying the conditions of Definition~\ref{Pgen} is greater than or equal to the number of partitions of $n$ with parts less than or equal to $3N$ satisfying the conditions of the definition of $G_{m,3N}(a,b,q)$.}

Conjecture~\ref{conjecture} calls for additional exploration. The authors suspect the existence of an injective proof of Conjecture~\ref{conjecture} in the spirit of \cite{BerkovichKeith}. We intend on addressing these conjectures in future work.

It has been mentioned in the introduction that Theorem~\ref{newcapparellicomp} is a corollary of Boulet's work \cite{Boulet}. Her work needs to be extended to deal with our refinement, Theorem~\ref{refined}. One of the interesting implications of this extension is:

\begin{theorem}\label{BouletConj} For $N$, $i$ and $j$ non-negative integers
\begin{equation*}\sum_{l=0}^{N} \sum_{k=0}^{l} (-1)^k \genfrac{[}{]}{0pt}{}{N}{l}_{q^2}\genfrac{[}{]}{0pt}{}{N-l}{i-k}_{q^2}\genfrac{[}{]}{0pt}{}{N-l}{j-k}_{q^2}\genfrac{[}{]}{0pt}{}{l}{k}_{q^2} q^{W(N,i,j,k,l)}= (-q;q)_{N-i-j} \genfrac{[}{]}{0pt}{}{N}{i,\ j}_{q^2}\end{equation*}
where $W(N,i,j,k,l):=3k^2+2(N-i-j-l)k +l$.
\end{theorem}
\noindent
Authors are going to address these new results in a separate paper. We leave the direct $q$-hypergeometric proof of Theorem~\ref{BouletConj} as an exercise to a motivated reader.

Finally, it is also of interest to understand the Affine Lie-theoretic connection of the new companion to Capparelli's theorem.

\section{Acknowledgement}

The authors would like to thank Krishnaswami Alladi, George E. Andrews, and Andrew V. Sills for their kind interest and encouragement. The authors also thank Frank Patane, Keith Grizzell, and Elizabeth Loew for their very helpful comments on the manuscript.

\end{document}